# CONVERGENCE OF MARKOV PROCESSES NEAR SADDLE FIXED POINTS

By Amanda G. Turner

*University of Cambridge*

We consider sequences $(X_t^N)_{t\geq 0}$ of Markov processes in two dimensions whose fluid limit is a stable solution of an ordinary differential equation of the form $\dot{x}_t = b(x_t)$, where $b(x) = \begin{pmatrix} -\mu & 0 \\ 0 & \lambda \end{pmatrix} x + \tau(x)$ for some $\lambda, \mu > 0$ and $\tau(x) = O(|x|^2)$. Here the processes are indexed so that the variance of the fluctuations of $X_t^N$ is inversely proportional to $N$. The simplest example arises from the OK Corral gunfight model which was formulated by Williams and McIlroy [*Bull. London Math. Soc.* **30** (1998) 166–170] and studied by Kingman [*Bull. London Math. Soc.* **31** (1999) 601–606]. These processes exhibit their most interesting behavior at times of order $\log N$ so it is necessary to establish a fluid limit that is valid for large times. We find that this limit is inherently random and obtain its distribution. Using this, it is possible to derive scaling limits for the points where these processes hit straight lines through the origin, and the minimum distance from the origin that they can attain. The power of $N$ that gives the appropriate scaling is surprising. For example if $T$ is the time that $X_t^N$ first hits one of the lines $y = x$ or $y = -x$, then

$$N^{\mu/(2(\lambda+\mu))}|X_T^N| \Rightarrow |Z|^{\mu/(\lambda+\mu)},$$

for some zero mean Gaussian random variable $Z$.

**1. Introduction.** The fluid limit theorem is a powerful result which shows that, under certain conditions, sequences of Markov processes converge to solutions of ordinary differential equations. We are interested in situations where the differential equation can be written in the form

$$\dot{x}_t = Bx_t + \tau(x_t), \tag{1}$$

for some matrix $B$, where $\tau(x) = O(|x|^2)$ is twice continuously differentiable. These differential equations have been studied extensively in the dynamical systems literature, with the aim of finding precise relationships between solutions of these differential equations and solutions of the corresponding linear differential equations

$$\dot{y}_t = By_t. \tag{2}$$

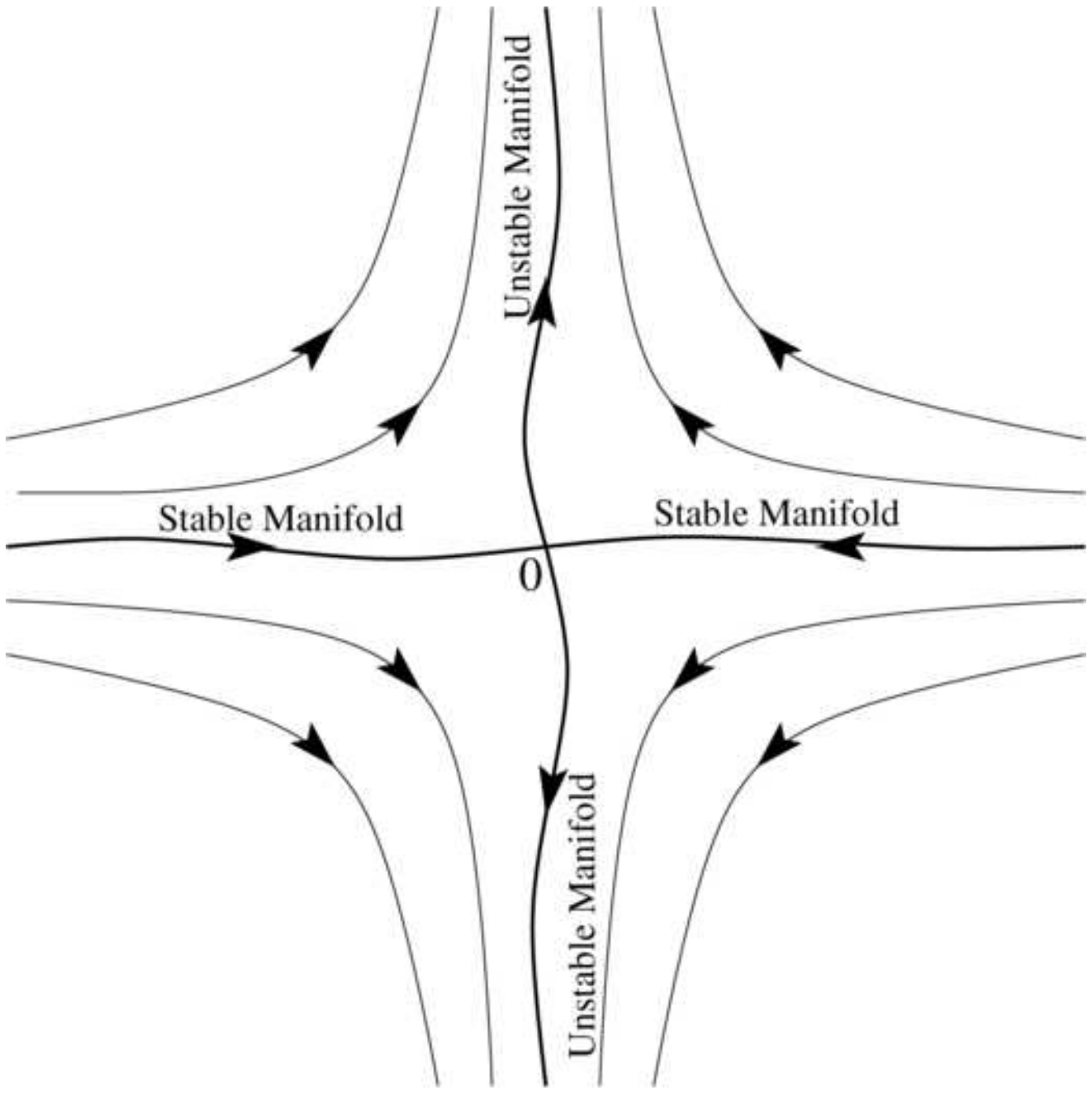

Fig. 1. *The phase portrait of an ordinary differential equation having a saddle fixed point at the origin.*

We restrict ourselves to the two dimensional case where the origin is a saddle fixed point of the system, that is, $B$ has eigenvalues $\lambda, -\mu$, with $\lambda, \mu > 0$. The phase portrait of (1) in the neighborhood of the origin is shown in Figure 1.

In particular, there exists some $x_0 \neq 0$ such that $\phi_t(x_0) \to 0$ as $t \to \infty$, where $\phi$ is the flow associated with the ordinary differential equation (1). The set of such $x_0$ is the stable manifold. There also exists some $x_\infty$ such that $\phi_t^{-1}(x_\infty) \to 0$ as $t \to \infty$. The set of such $x_\infty$ is the unstable manifold. The saddle case is interesting in this setting as it is the only case in two dimensions where there is both a stable and an unstable manifold.

Fix an $x_0$ in the stable manifold and consider sequences of Markov processes with initial condition $X_0^N = x_0$, where the processes are indexed so that the variance of the fluctuations of $X_t^N$ is inversely proportional to $N$. The fluid limit theorem tells us that for fixed values of $t$, $X_t^N \to \phi_t(x_0)$ as $N \to \infty$. However, if we allow the value of $t$ to grow with $N$ as $N \to \infty$, we shall see that $X_t^N$ deviates from the stable solution to a limit which is inherently random, before converging to an unstable solution (see Figure 2).

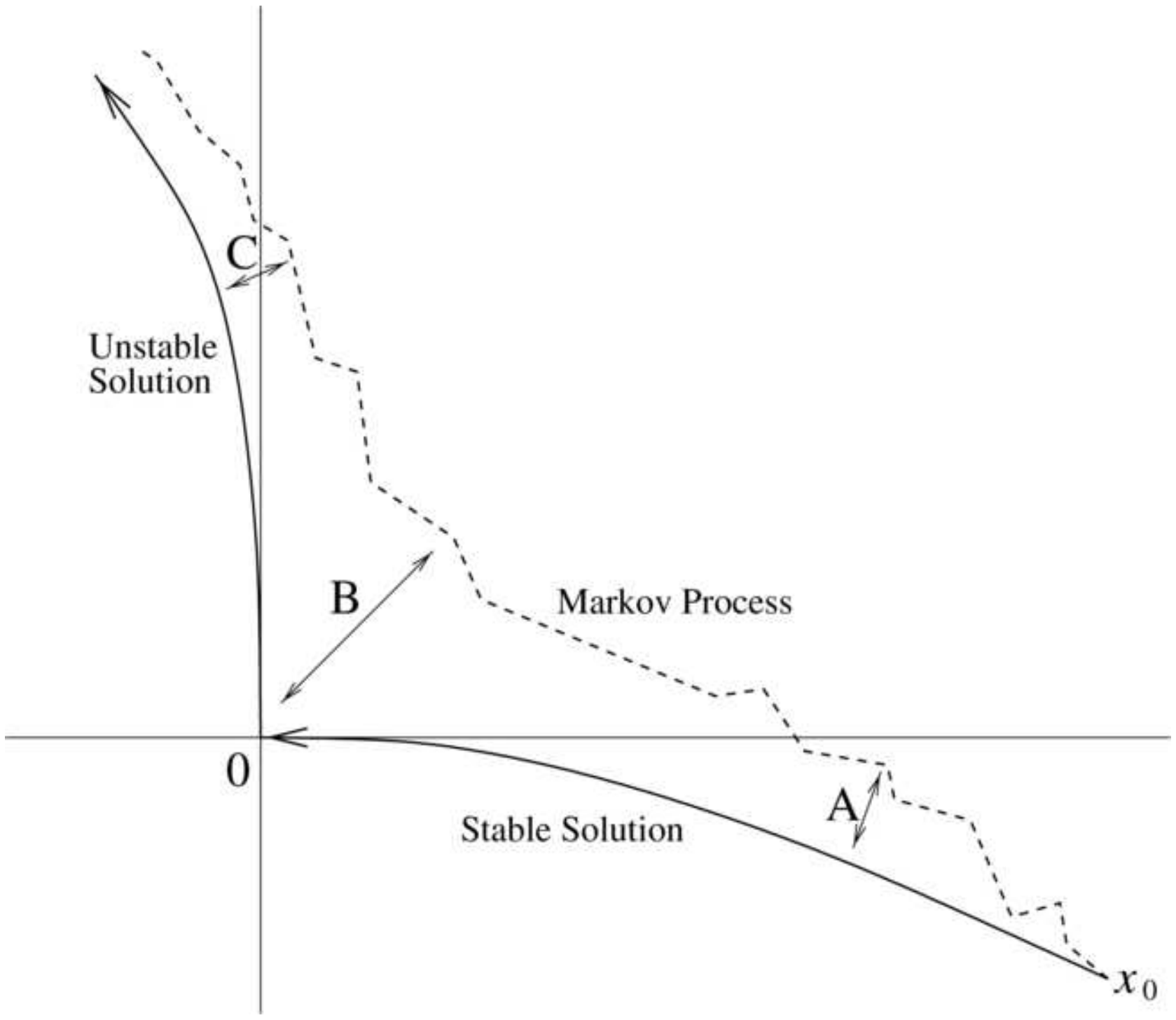

FIG. 2. *Diagram showing how the Markov process $X_t^N$ deviates from the stable solution $\phi_t(x_0)$ for large values of $t$.*

More precisely, we observe three different types of behavior depending on the time scale:

(A) On compact time intervals, $X_t^N$ converges to the stable solution of (1), the fluctuations around this limit being of order $N^{-1/2}$.

(B) There exists some $\bar{x}_0 \neq 0$, depending only on $x_0$, and a Gaussian random variable $Z_\infty$ such that if $t$ lies in the interval $[R, \frac{1}{2\lambda} \log N - R]$, then

$$X_t^N = \bar{x}_0 e^{-\mu t}(e_1 + \varepsilon_1) + N^{-1/2} Z_\infty e^{\lambda t}(e_2 + \varepsilon_2)$$

for some $\varepsilon_i(t, N) \to 0$ uniformly in $t$ in probability as $R, N \to \infty$, where $\{e_1, e_2\}$ is the standard basis for $\mathbb{R}^2$. In other words, $X_t^N$ can be approximated by the solution to the linear ordinary differential equation (2) starting from the random point $\binom{\bar{x}_0}{N^{-1/2} Z_\infty}$.

(C) On time intervals of a fixed length around $\frac{1}{2\lambda} \log N$, $X_t^N$ converges to the unstable solution of (1).

The most interesting behavior occurs on time intervals of fixed lengths around $\frac{1}{2(\lambda+\mu)}\log N$, as for these values of $t$ the two terms $\bar{x}_0 e^{-\mu t}$ and $N^{-1/2}Z_\infty e^{\lambda t}$ are of the same order. By considering

$$\bar{x}_0 e^{-\mu t} e_1 + N^{-1/2} Z_\infty e^{\lambda t} e_2,$$

we show in Section 7 that it is at these times that $X_t^N$ crosses all the straight lines passing through 0, and also that $|X_t^N|$ attains its minimum value when $t$ is in this range. The distance from the origin of $X_t^N$ for these values of $t$ is of order $N^{-\mu/(2(\lambda+\mu))}$, which gives us surprising scaling limits for the points at which $X_t^N$ intersects various straight lines, and for $\inf |X_t^N|$.

In order to study the Markov processes at times of order $\log N$, it is necessary to establish a strong form of the fluid limit theorem that is valid for large times. The key idea is to show that for $N$ and $t_0$ sufficiently large, the process $(X_t^N)_{t\geq t_0}$ is close to $(\phi_{t-t_0}(X_{t_0}^N))_{t\geq t_0}$. This is done in Section 2 in the case when (1) is linear and $X_t^N$ is a pure jump Markov process, in Section 5 for general pure jump Markov processes, and in Section 6 for continuous diffusion processes. In Sections 3 and 4 we look at the process $(\phi_{t-t_0}(X_{t_0}^N))_{t\geq t_0}$ for large values of $N$ and $t_0$, which then enables us to obtain scaling limits for the process $X_t^N$. The same idea can be used to obtain fluid limit theorems for arbitrary matrices $B$, for example, with eigenvalues having the same sign, or in higher dimensions, however an analysis of the solutions of the underlying differential equation is required, which we do not go into here.

The simplest example of this type of behavior arises from the OK Corral gunfight model which was formulated by Williams and McIlroy [7] and studied by Kingman [5] and Kingman and Volkov [6]. Two lines of gunmen face each other, there initially being $N$ on each side. Each gunman fires lethal gunshots at times of a Poisson process with rate 1 until either there is no one left on the other side or he is killed. The process terminates when all the gunmen on one side are dead. It is shown by Kingman that if $S^N$ is the number of survivors when the process terminates, then

$$N^{-3/4} S^N \Rightarrow 2^{3/4} |Z|^{1/2},$$

where $Z \sim N(0, \frac{1}{3})$. It is the occurrence of the unexpected power of $N$ that interested the above authors in the problem. By using our scaling limits we rederive this result in Section 2.1 and show that it is a special case of a much more general phenomenon, and that in fact by a suitable choice of $B$, every number in the interval $(\frac{1}{2}, 1)$ may be obtained as a power of $N$ in this way. An application of the nonlinear case to a model of two competing species is given in Section 7.

**2. The linear case.** In this section we restrict ourselves to sequences of Markov processes in the special case where equation (1) is linear. We begin by describing the conditions under which a limit theorem exists for large times and then establish the exact limit by means of an appropriate martingale inequality. In Section 2.1 this result is used to derive scaling limits for the points where these processes hit straight lines through the origin and we use this to obtain a solution to the OK Corral problem.

The fluid limit theorem that we state below is widely known and has been the object of many works. We use the formulation found in [2].

Let $(X_t^N)_{t\geq 0}$ be a sequence of pure jump Markov processes, starting from $x_0$ and taking values in some subsets $I^N$ of $\mathbb{R}^2$, with Lévy kernels $K^N(x,dy)$. Let $S$ be an open subset of $\mathbb{R}^2$ with $x_0 \in S$, and set $S^N = I^N \cap S$. For $x \in S^N$ and $\theta \in (\mathbb{R}^2)^*$, define the Laplace transform corresponding to Lévy kernel $K^N(x,dy)$ by

$$m^N(x,\theta) = \int_{\mathbb{R}^2} e^{\langle \theta, y\rangle} K^N(x,dy).$$

We assume that there is a limit kernel $K(x,dy)$ defined for $x \in S$ with corresponding Laplace transform $m(x,\theta)$ with the following properties:

(a) There exists a constant $\eta_0 > 0$ such that $m(x,\theta)$ is uniformly bounded for all $x \in S$ and $|\theta| \leq \eta_0$.

(b) As $N \to \infty$,

$$\sup_{x\in S^N} \sup_{|\theta|\leq \eta_0} \left| \frac{m^N(x, N\theta)}{N} - m(x,\theta) \right| \to 0.$$

Set $b(x) = m'(x,0)$ where $'$ denotes differentiation in $\theta$. Suppose that $b$ is Lipschitz on $S$ so that $b$ has an extension to a Lipschitz vector field $\tilde{b}$ on $\mathbb{R}^2$. Then there is a unique solution $(x_t)_{t\geq 0}$ to the ordinary differential equation $\dot{x}_t = \tilde{b}(x_t)$ starting from $x_0$. Suppose that $S$ contains a neighborhood of the path $(x_t)_{t\geq 0}$. By stopping $X_t^N$ at the first time it leaves $S$, if necessary, we may assume that $X_t^N$ remains in $\overline{S}$ for all $t \geq 0$. Under these assumptions, for all $t_0 \geq 0$ and $\delta > 0$,

$$\limsup_{N\to\infty} N^{-1} \log \mathbb{P}\left( \sup_{t\leq t_0} |X_t^N - x_t| \geq \delta \right) < 0.$$

Suppose additionally that:

(c) $b$ is $C^1$ on $S$ and

$$\sup_{x\in S^N} N^{1/2} |b^N(x) - b(x)| \to 0,$$

where $b^N(x) = m^{N\prime}(x,0)$.

(d) $a$, defined by $a(x) = m''(x, 0)$, is Lipschitz on $S$.

It follows from the above that for any $\eta < \eta_0$ there exists a constant $A$ such that

$$\sup_{x \in S^N} \sup_{|\theta| \leq \eta} N|m^{N''}(x, N\theta)| \leq A. \tag{3}$$

Let $\gamma_t^N = N^{1/2}(X_t^N - x_t)$. Then for any $t \geq 0$, $\gamma_t^N \Rightarrow \gamma_t$ as $N \to \infty$, where $(\gamma_t)_{t \geq 0}$ is the unique solution to the linear stochastic differential equation

$$d\gamma_t = \sigma(x_t)\, dW_t + \nabla b(x_t)\gamma_t\, dt \tag{4}$$

starting from 0, $W$ a Brownian motion in $\mathbb{R}^2$, and $\sigma \in \mathbb{R}^2 \otimes (\mathbb{R}^2)^*$ satisfying $\sigma(x)\sigma(x)^* = a(x)$. The distribution of $(\gamma_t)_{t \geq 0}$ does not depend on the choice of $\sigma$.

We are interested in the case where $b(x) = Bx$ for some matrix $B = \left(\begin{smallmatrix} -\mu & 0 \\ 0 & \lambda \end{smallmatrix}\right)$, $\mu, \lambda > 0$.

Let $\phi_t(x)$ be the solution to the ordinary differential equation

$$\dot{\phi}_t(x) = b(\phi_t(x)), \qquad \phi_0(x) = x. \tag{5}$$

In the linear case we can solve (5) explicitly to get $\phi_t(x) = e^{Bt}x$. We concentrate on processes where the initial condition is chosen to be $x_0 = (x_{0,1}, 0)$ with $x_{0,1} \neq 0$, so that $x_t = \phi_t(x_0) \to 0$ as $t \to \infty$. We shall show that for sufficiently large values of $N$ and $t_0$, $X_t^N$ is in some sense close to $\phi_{t-t_0}(X_{t_0}^N)$ for $t \geq t_0$.

Introduce random measures $\mu^N$ and $\nu^N$ on $(0, \infty) \times \mathbb{R}^2$, given by

$$\mu^N = \sum_{\Delta X_t^N \neq 0} \delta_{(t, \Delta X_t^N)},$$

$$\nu^N(dt, dy) = K^N(X_{t-}^N, dy)\, dt,$$

where $\delta_{(t,y)}$ denotes the unit mass at $(t, y)$ and $\Delta X_t^N = X_t^N - X_{t-}^N$.

Let $f(t, x) = e^{-Bt}(x - \phi_{t-t_0}(X_{t_0}^N))$, for $t \geq t_0$. By Itô's formula,

$$f(t, X_t^N) = f(t_0, X_{t_0}^N) + M_t^{B,N} - M_{t_0}^{B,N} + \int_{t_0}^t \left(\frac{\partial f}{\partial t} + K^N f\right)(s, X_{s-}^N)\, ds,$$

where

$$\frac{\partial f}{\partial t} = -Be^{-Bt}x,$$

$$\begin{aligned} K^N f(s, x) &= \int_{\mathbb{R}^2} (f(s, x + y) - f(s, x))K^N(x, dy) \\ &= \int_{\mathbb{R}^2} e^{-Bs} y K^N(x, dy) \\ &= e^{-Bs} b^N(x) \end{aligned}$$

and

$$M_t^{B,N} = \int_{(0,t]\times\mathbb{R}^2} (f(s, X_{s-}^N + y) - f(s, X_{s-}^N))(\mu^N - \nu^N)(ds, dy)$$
$$= \int_{(0,t]\times\mathbb{R}^2} e^{-Bs} y(\mu^N - \nu^N)(ds, dy).$$

So if $t \geq t_0$, then

$$\begin{aligned} & e^{-Bt}(X_t^N - \phi_{t-t_0}(X_{t_0}^N)) \\ & \quad = M_t^{B,N} - M_{t_0}^{B,N} + \int_{t_0}^t e^{-Bs}(b^N(X_{s-}^N) - b(X_{s-}^N))\,ds. \end{aligned} \tag{6}$$

Lemma 2.1. *There exists some constant $C$ such that*

$$\mathbb{E}\Big(\sup_{t\geq t_0} e^{-\lambda t}|e^{Bt}(M_t^{B,N} - M_{t_0}^{B,N})|\Big) \leq CN^{-1/2}e^{-\lambda t_0}.$$

Proof. By the product rule,

$$\begin{aligned} & e^{(B-\lambda I)t}(M_t^{B,N} - M_{t_0}^{B,N}) \\ & \quad = \int_{t_0}^t (B - \lambda I)e^{(B-\lambda I)s}(M_s^{B,N} - M_{t_0}^{B,N})\,ds \\ & \qquad + \int_{t_0}^t \int_{\mathbb{R}^2} e^{-\lambda s} y(\mu^N - \nu^N)(dy, ds) \end{aligned}$$

and hence

$$\begin{aligned} & \mathbb{E}\Big(\sup_{t\geq t_0} e^{-\lambda t}|e^{Bt}(M_t^{B,N} - M_{t_0}^{B,N})|\Big) \\ & \quad \leq \mathbb{E}\Big(\sup_{t\geq t_0} \int_{t_0}^t (\lambda + \mu)e^{-(\lambda+\mu)s}|(M_s^{B,N} - M_{t_0}^{B,N})_1|\,ds\Big) \\ & \qquad + \mathbb{E}\Big(\sup_{t\geq t_0}\Big|\int_{t_0}^t \int_{\mathbb{R}^2} e^{-\lambda s} y(\mu^N - \nu^N)(dy, ds)\Big|\Big) \\ & \quad \leq \int_{t_0}^\infty (\lambda + \mu)e^{-(\lambda+\mu)s}(\mathbb{E}(M_s^{B,N} - M_{t_0}^{B,N})_1^2)^{1/2}\,ds \\ & \qquad + \mathbb{E}\Big(\sup_{t\geq t_0}\Big|\int_{t_0}^t \int_{\mathbb{R}^2} e^{-\lambda s} y(\mu^N - \nu^N)(dy, ds)\Big|^2\Big)^{1/2}. \end{aligned}$$

Since

$$\mathbb{E}\int_0^t \int_{\mathbb{R}^2} |e^{-\lambda s} y|\nu^N(ds, dy) < \infty$$

for all $t \geq 0$, the process

$$\left(\int_0^t \int_{\mathbb{R}^2} e^{-\lambda s} y(\mu^N - \nu^N)(dy, ds)\right)_{t \geq 0}$$

is a martingale, and hence, by Doob's $L^2$ inequality

$$\mathbb{E}\left(\sup_{t \geq t_0} \left| \int_{t_0}^t \int_{\mathbb{R}^2} e^{-\lambda s} y(\mu^N - \nu^N)(dy, ds)\right|^2\right)$$
$$\leq 4 \sup_{t \geq t_0} \mathbb{E}\left(\left| \int_{t_0}^t \int_{\mathbb{R}^2} e^{-\lambda s} y(\mu^N - \nu^N)(dy, ds)\right|^2\right).$$

Now

$$\begin{aligned}\mathbb{E}((M_t^{B,N} - M_{t_0}^{B,N})_1^2) &= \mathbb{E} \int_{t_0}^t \int_{\mathbb{R}^2} e^{2\mu s} y_1^2 \nu^N(dy, ds) \\ &\leq \mathbb{E} \int_{t_0}^t e^{2\mu s} |m^{N\prime\prime}(X_{s-}^N, 0)|\, ds \\ &\leq \frac{e^{2\mu t} A}{2\mu N},\end{aligned}$$

where $A$ is defined in (3). Similarly

$$\mathbb{E}\left(\left| \int_{t_0}^t \int_{\mathbb{R}^2} e^{-\lambda s} y(\mu^N - \nu^N)(dy, ds)\right|^2\right) \leq \frac{e^{-2\lambda t_0} A}{2\lambda N}.$$

Hence

$$\begin{aligned}&\mathbb{E}\left(\sup_{t \geq t_0} e^{-\lambda t} |e^{Bt}(M_t^{B,N} - M_{t_0}^{B,N})|\right) \\ &\quad \leq \int_{t_0}^\infty (\lambda + \mu) e^{-\lambda s} \left(\frac{A}{2\mu N}\right)^{1/2} ds + e^{-\lambda t_0} \left(\frac{2A}{\lambda N}\right)^{1/2} \\ &\quad \leq \frac{A^{1/2}(\lambda + \mu + 2(\lambda\mu)^{1/2})}{\lambda (2\mu)^{1/2}} N^{-1/2} e^{-\lambda t_0}.\end{aligned}$$

□

THEOREM 2.2. *For all $\varepsilon > 0$,*

$$\lim_{t_0 \to \infty} \limsup_{N \to \infty} \mathbb{P}\left(\sup_{t \geq t_0} e^{-\lambda t} |X_t^N - \phi_{t - t_0}(X_{t_0}^N)| > N^{-1/2} \varepsilon\right) = 0.$$

PROOF. Let $N_0$ be sufficiently large that $\sup_{N \geq N_0} N^{1/2} \|b^N - b\| < \lambda\varepsilon/2$ and set

$$\Omega_{N,t_0} = \left\{\sup_{t \geq t_0} |e^{-\lambda t} e^{Bt}(M_t^{B,N} - M_{t_0}^{B,N})| \leq N^{-1/2} \frac{\varepsilon}{2}\right\}.$$

By (6), on the set $\Omega_{N,t_0}$ with $N \geq N_0$,

$$\begin{aligned}\sup_{t\geq t_0} e^{-\lambda t}|X_t^N - \phi_{t-t_0}(X_{t_0}^N)| &\leq \sup_{t\geq t_0} |e^{-\lambda t}e^{Bt}(M_t^{B,N} - M_{t_0}^{B,N})| \\ &\quad + \sup_{t\geq t_0} e^{-\lambda t}\int_{t_0}^t |e^{B(t-s)}|\|b^N - b\|\,ds \\ &\leq N^{-1/2}\varepsilon.\end{aligned}$$

Hence

$$\begin{aligned}\limsup_{N\to\infty}\mathbb{P}\Big(\sup_{t\geq t_0} e^{-\lambda t}|X_t^N - \phi_{t-t_0}(X_{t_0}^N)| > N^{-1/2}\varepsilon\Big) &\leq \limsup_{N\to\infty}\mathbb{P}(\Omega_{N,t_0}^c) \\ &\leq \frac{2Ce^{-\lambda t_0}}{\varepsilon} \\ &\to 0\end{aligned}$$

as $t_0 \to \infty$, where the second inequality follows by Markov's inequality and Lemma 2.1. □

Let $Z_\infty \sim N(0, \sigma_\infty^2)$, where

$$\sigma_\infty^2 = \int_0^\infty e^{-2\lambda s} a(x_s)_{2,2}\,ds.$$

THEOREM 2.3. *The following converge in probability as $N \to \infty$.*

(i)

$$\sup_{t\leq t_N} |e^{\mu t}X_{t,1}^N - x_{0,1}| \to 0$$

*for any sequence $t_N \to \infty$ with $e^{(\lambda+\mu)t_N} = O(N^{1/2})$;*

(ii)

$$\sup_{t\geq t_N} N^{1/2}e^{-\lambda t}|X_{t,1}^N| \to 0$$

*for any sequence $t_N$ with $e^{(\lambda+\mu)t_N} = \omega(N^{1/2})$;*

(iii)

$$\sup_{t_1,t_2\geq t_N} N^{1/2}|e^{-\lambda t_1}X_{t_1,2}^N - e^{-\lambda t_2}X_{t_2,2}^N| \to 0$$

*for any sequence $t_N \to \infty$.*

*Furthermore, if $\sigma_\infty \neq 0$, then*

$$N^{1/2}e^{-\lambda t}X_{t,2}^N \Rightarrow Z_\infty$$

*as $t, N \to \infty$.*

REMARK 2.4. Given any sequence of times $t_N \to \infty$ as $N \to \infty$, by the Skorohod representation theorem, it is possible to choose a sample space in which $Z_\infty^N = N^{1/2} e^{-\lambda t_N} X_{t_N,2}^N \to Z_\infty$ almost surely as $N \to \infty$. In this case the above result can be expressed as

$$X_t^N = x_{0,1} e^{-\mu t}(e_1 + \varepsilon_1) + N^{-1/2} Z_\infty e^{\lambda t}(e_2 + \varepsilon_2) \tag{7}$$

where $\varepsilon_i = \varepsilon_i(N,t) \to 0$, uniformly in t, in probability as $N \to \infty$.

PROOF OF THEOREM 2.3. For any fixed $t_0$, $\sup_{t \le t_0} |e^{\mu t} X_{t,1}^N - x_{0,1}| \to 0$ as an immediate consequence of the fluid limit theorem. For (i) it is therefore sufficient to show that for any $\varepsilon > 0$, $\lim_{t_0 \to \infty} \limsup_{N \to \infty} \mathbb{P}(\sup_{t_0 \le t_N} |e^{\mu t} X_{t,1}^N - x_{0,1}| > \varepsilon) = 0$. Now if $t \ge t_0$, $\phi_{t-t_0}(X_{t_0}^N) = e^{B(t-t_0)} X_{t_0}^N = e^{B(t-t_0)}(x_{t_0} + N^{-1/2} \gamma_{t_0}^N)$. Since $x_0 = x_{0,1} e_1$, $e^{B(t-t_0)} x_{t_0} = e^{-\mu t} x_0$. Hence

$$\phi_{t-t_0}(X_{t_0}^N) = e^{-\mu t}(x_0 + N^{-1/2} e^{\mu t_0} \gamma_{t_0,1}^N e_1) + N^{-1/2} e^{\lambda t} e^{-\lambda t_0} \gamma_{t_0,2}^N e_2,$$

and so

$$X_{t,1}^N = e^{-\mu t}(x_{0,1} + N^{-1/2} e^{\mu t_0} \gamma_{t_0,1}^N + e^{(\lambda+\mu)t} e^{-\lambda t}(X_t^N - \phi_{t-t_0}(X_{t_0}^N))_1)$$

and

$$X_{t,2}^N = N^{-1/2} e^{\lambda t}(e^{-\lambda t_0} \gamma_{t_0,2}^N + N^{1/2} e^{-\lambda t}(X_t^N - \phi_{t-t_0}(X_{t_0}^N))_2). \tag{8}$$

Let $N \to \infty$ and then $t_0 \to \infty$. Statements (i)–(iii) follow by Theorem 2.2 and the fact that $\gamma_{t_0}^N \Rightarrow \gamma_{t_0}$, a Gaussian random variable.

For the last part, note that by (4),

$$e^{-\lambda t_0} \gamma_{t_0,2}^N \Rightarrow e^{-\lambda t_0} \gamma_{t_0,2} = \int_0^{t_0} e^{-\lambda s} \langle e_2, \sigma(x_s)\, dW_s \rangle$$

as $N \to \infty$. Since

$$\begin{aligned} \int_0^\infty |e^{-\lambda s} e_2^* \sigma(x_s)|^2\, ds &\le \int_0^\infty e^{-2\lambda s} |a(x_s)|\, ds \\ &\le \frac{A}{2\lambda}, \end{aligned}$$

where $e_i^*$ is the transpose of $e_i$, $e^{-\lambda t_0} \gamma_{t_0,2} \to Z_\infty$ almost surely as $t_0 \to \infty$, for

$$Z_\infty = \left( \int_0^\infty e^{-\lambda t} \sigma(x_t)\, dW_t \right)_2 \sim N(0, \sigma_\infty^2).$$

The result follows by (8) and Theorem 2.2. □

2.1. *Applications.* Applications will be dealt with more fully in Section 7. However, we illustrate here how the above result can be used to study the first time that $X_t^N$ hits $l_\theta$ or $l_{-\theta}$, the straight lines passing through the origin at angles $\theta$ and $-\theta$, where $\theta \in (0, \frac{\pi}{2})$, as $N \to \infty$. As $X_t^N$ is not continuous, we define the time that $X_t^N$ first intersects one of the lines $l_{\pm\theta}$ as

$$T_\theta^N = \inf\left\{t \ge 0 : \left|\frac{X_{t-,2}^N}{X_{t-,1}^N}\right| \le |\tan\theta| \text{ and } \left|\frac{X_{t,2}^N}{X_{t,1}^N}\right| \ge |\tan\theta|\right\}.$$

Let

$$t_N = \frac{1}{2(\lambda+\mu)} \log N$$

and

$$c_\theta = \frac{1}{\lambda+\mu} \log\left|\frac{x_{0,1}\tan\theta}{Z_\infty}\right|.$$

THEOREM 2.5. *Under the assumptions of Theorem 2.3*

$$T_\theta^N - t_N \Rightarrow c_\theta \tag{9}$$

*and*

$$N^{\mu/(2(\lambda+\mu))}|X_{T_\theta^N}^N| \Rightarrow |\sec\theta||\tan\theta|^{-\mu/(\lambda+\mu)}|x_0|^{\lambda/(\lambda+\mu)}|Z_\infty|^{\mu/(\lambda+\mu)} \tag{10}$$

*as* $N \to \infty$.

PROOF. For simplicity, we work in a sample space in which $Z_\infty^N \to Z_\infty$ almost surely. Define $\varepsilon_i$ as in Remark 2.4. By observing that

$$x_{0,1}e^{-\mu t}e_1 + N^{-1/2}Z_\infty e^{\lambda t}e_2$$

first intersects one of the lines $l_{\pm\theta}$ at time $t = t_N + c_\theta$, given any $\varepsilon > 0$,

$$\begin{aligned}
&\mathbb{P}(T_\theta^N \le t_N + c_\theta - \varepsilon) \\
&\quad \le \mathbb{P}\left(\sup_{t \le t_N + c_\theta - \varepsilon}\left|\frac{X_{t,2}^N}{X_{t,1}^N}\right| \ge |\tan\theta|\right) \\
&\quad = \mathbb{P}\left(\sup_{t \le t_N + c_\theta - \varepsilon}\left|\frac{x_{0,1}e^{-\mu t}\varepsilon_{1,2} + N^{-1/2}Z_\infty e^{\lambda t}(1+\varepsilon_{2,2})}{x_{0,1}e^{-\mu t}(1+\varepsilon_{1,1}) + N^{-1/2}Z_\infty e^{\lambda t}\varepsilon_{2,1}}\right| \ge |\tan\theta|\right) \\
&\quad \to 0
\end{aligned}$$

as $N \to \infty$. Similarly,

$$\mathbb{P}(T_\theta^N \ge t_N + c_\theta + \varepsilon) \le \mathbb{P}\left(\inf_{t \ge t_N + c_\theta + \varepsilon}\left|\frac{X_{t,2}^N}{X_{t,1}^N}\right| \le |\tan\theta|\right) \to 0.$$

The result follows immediately. $\square$

REMARK 2.6. The sign of $Z_\infty$ determines whether $X_t^N$ hits $l_\theta$ or $l_{-\theta}$ at time $T_\theta^N$. Since $Z_\infty$ is a Gaussian random variable with mean 0, each event occurs with probability $\frac{1}{2}$.

EXAMPLE 2.7 (The OK Corral problem). The OK Corral process is a $\mathbb{Z}^2$ valued process $(U_t^N, V_t^N)$ used to model the famous gunfight where $U_t^N$ and $V_t^N$ are the number of gunmen on each side and $U_0^N = V_0^N = N$. Each gunman fires lethal gunshots at times of a Poisson process with rate 1 until either there is no one left on the other side or he is killed. The transition rates are

$$(u,v) \to \begin{cases} (u-1,v), & \text{at rate } v, \\ (u,v-1), & \text{at rate } u \end{cases}$$

until $uv = 0$.

The process terminates when all the gunmen on one side are dead. We are interested in the number of gunmen surviving when the process terminates, for large values of $N$.

This model was formulated by Williams and McIlroy [7] and later studied by Kingman [5] and subsequently Kingman and Volkov [6].

Let $X_t^N = (U_t^N, V_t^N)/N$. This gives a sequence of pure jump Markov processes, starting from $x_0 = (1,1)$, with Lévy kernels

$$K^N(x,dy) = Nx_2\delta_{(-1/N,0)} + Nx_1\delta_{(0,-1/N)}.$$

If we let

$$K(x,dy) = x_2\delta_{(-1,0)} + x_1\delta_{(0,-1)},$$

then

$$m(x,\theta) = x_2e^{-\theta_1} + x_1e^{-\theta_2} = \frac{m^N(x,N\theta)}{N},$$

$$b(x) = \begin{pmatrix} 0 & -1 \\ -1 & 0 \end{pmatrix} x = b^N(x)$$

and

$$a(x) = \begin{pmatrix} x_2 & 0 \\ 0 & x_1 \end{pmatrix}.$$

So, under a rotation by $\frac{\pi}{4}$, the conditions required for Theorem 2.5 are satisfied, with $\lambda = \mu = 1$. In the original coordinates, the process terminates when $X_t^N$ hits the $x$ or $y$ axes. Under the rotation, this corresponds to hitting $l_{\pm\frac{\pi}{4}}$. Hence, if the OK Corral process terminates at time $T_N$ and there are $S^N$ survivors, then $T_N = T_{\pi/4}^N$ and $S^N = N|X^N_{T^N_{\pi/4}}|$, and so

$$T_N - \tfrac{1}{4}\log N \Rightarrow \tfrac{1}{4}\log 2 - \tfrac{1}{2}\log|Z_\infty|$$

and

$$N^{-3/4}S^N \Rightarrow 2^{3/4}|Z_\infty|^{1/2},$$

where $Z_\infty \sim N(0, \frac{1}{3})$. The limiting distribution of $N^{-3/4}S^N$ is the one obtained by Kingman in [5].

REMARK 2.8. It is remarked by Kingman [5] that it is the occurrence of the surprising power of $N$ that makes the OK Corral process of interest. Theorem 2.5 shows that this is a special case of a more general phenomenon, and in fact, by a suitable choice of $\frac{\lambda}{\mu}$, every number in the interval $(\frac{1}{2}, 1)$ may be obtained as a power of $N$ in this way.

**3. Linearization of the limit process.** We now turn to the general case where $b(x) = Bx + \tau(x)$ for $B = \left(\begin{smallmatrix} -\mu & 0 \\ 0 & \lambda \end{smallmatrix}\right)$, $\mu, \lambda > 0$, and $\tau : \mathbb{R}^2 \to \mathbb{R}^2$, twice continuously differentiable, with $\tau(0) = \nabla\tau(0) = 0$. Let $\phi_t(x)$ be the solution to the ordinary differential equation

$$\dot{\phi}_t(x) = b(\phi_t(x)), \qquad \phi_0(x) = x. \tag{11}$$

This section consists of a technical calculation which expresses $\phi_t(x)$ in a linear form.

We are interested in the behavior of solutions starting near the stable manifold. Lemma 3.2 proves the existence of the stable manifold and establishes the limiting behavior of a stable solution. First order behavior is investigated in Lemma 3.3, and these results are then used in Theorem 3.4 to express solutions near the stable manifold in the required linear form. Theorem 3.5 shows that over large time periods, solutions starting near the stable manifold approach the unstable manifold.

Throughout this section we use the following classical planar linearization theorem due to Hartman [4].

THEOREM 3.1. *There exists a $C^1$ diffeomorphism $h : U \to V = h(U)$, defined on an open neighborhood $U$ of the origin, with uniformly Hölder continuous partial derivatives and having the form $h(x) = x + o(x)$ such that*

$$h(\phi_t(x)) = e^{Bt}h(x)$$

*for all $(t, x)$ with $\phi_t(x) \in U$.*

Pick $0 < \delta < 1$ sufficiently small that the ball of radius $\delta$ centered at the origin is contained in $U \cap V$. Since $h^{-1}(x) = x + o(x)$, and $\nabla h(x) = I + o(1)$ we can further ensure that $\delta$ is sufficiently small that

$$\sup_{0<|x|<\delta} (|h(x)/x| \vee |h^{-1}(x)/x|) < 2$$

and

$$\sup_{|x|<\delta}(|\nabla h(x) - I| \vee |\nabla h^{-1}(x) - I|) < 1/2.$$

LEMMA 3.2. *There exists an $x_0$ with $0 < |x_0| < \delta/8$ such that $\phi_t(x_0) \to 0$ as $t \to \infty$. Furthermore, for any such $x_0$, there exists some $\bar{x}_0$ with $0 < |\bar{x}_0| < \delta/4$ such that*

$$e^{\mu t}\phi_t(x_0) \to \begin{pmatrix} \bar{x}_0 \\ 0 \end{pmatrix}$$

*as $t \to \infty$, and*

$$|\phi_t(x_0)| \le 2|\bar{x}_0|e^{-\mu t} < \delta e^{-\mu t}/2$$

*for all $t \ge 0$.*

PROOF. Pick some $\bar{x}_0 \in \mathbb{R}$ with $0 < |\bar{x}_0| < \delta/16$ and define $x_0 = h^{-1}(\bar{x}_0, 0)$. Then

$$0 < |x_0| \le \sup_{0<|x|<\delta} |h^{-1}(x)/x||\bar{x}_0| < \frac{\delta}{8}$$

and

$$\phi_t(x_0) = h^{-1}\left(e^{Bt}\begin{pmatrix} \bar{x}_0 \\ 0 \end{pmatrix}\right) = h^{-1}\begin{pmatrix} e^{-\mu t}\bar{x}_0 \\ 0 \end{pmatrix} \to 0$$

as $t \to \infty$.

Conversely, given $x_0$ satisfying the above conditions, define $\bar{x}_0 = h(x_0)_1$. Note that because of the form of $h(x)$, $\bar{x}_0$ has the same sign as $x_{0,1}$. Since $e^{Bt}h(x_0) = h(\phi_t(x_0)) \to 0$ as $t \to \infty$, $h(x_0)_2 = 0$, and so

$$0 < |\bar{x}_0| = |h(x_0)| \le 2|x_0| < \delta/4.$$

Also

$$e^{\mu t}\phi_t(x_0) = e^{\mu t}h^{-1}\left(e^{Bt}\begin{pmatrix} \bar{x}_0 \\ 0 \end{pmatrix}\right) = e^{\mu t}\left(\begin{pmatrix} e^{-\mu t}\bar{x}_0 \\ 0 \end{pmatrix} + o(e^{-\mu t}\bar{x}_0)\right) \to \begin{pmatrix} \bar{x}_0 \\ 0 \end{pmatrix}$$

as $t \to \infty$, and

$$|\phi_t(x_0)| = \left|h^{-1}\left(e^{Bt}\begin{pmatrix} \bar{x}_0 \\ 0 \end{pmatrix}\right)\right| = \left|h^{-1}\begin{pmatrix} e^{-\mu t}\bar{x}_0 \\ 0 \end{pmatrix}\right| \le 2|\bar{x}_0|e^{-\mu t} < \frac{\delta}{2}e^{-\mu t}$$

for all $t \ge 0$. □

LEMMA 3.3. (i) *There exists some $D_0 \in (\mathbb{R}^2)^* \setminus \{\underline{0}\}$, where $\underline{0} = (0 \quad 0)$, such that*

$$e^{-\lambda t}\nabla\phi_t(x_0) \to \begin{pmatrix} \underline{0} \\ D_0 \end{pmatrix}$$

*as $t \to \infty$.*

(ii) *If* $|x| < \delta$ *and* $|\phi_t(x)| < \delta/2$*, then* $|\nabla\phi_t(x)| < 4e^{\lambda t}$.

(iii) *If* $|x| + |y| < \delta$ *and* $\sup_{0\leq\theta\leq 1} |\phi_t(x+\theta y)| < \delta/2$*, then there exist constants* $K \in \mathbb{R}$ *and* $0 < \alpha \leq 1$ *such that*

$$|\nabla\phi_t(x+y) - \nabla\phi_t(x)| \leq K e^{\lambda t(1+\alpha)} |y|^\alpha.$$

PROOF. (i) Let $D_0 = \nabla h_2(x_0) \in (\mathbb{R}^2)^* \setminus \{\underline{0}\}$. Then

$$\begin{aligned} e^{-\lambda t}\nabla\phi_t(x_0) &= \nabla h^{-1}\begin{pmatrix} e^{-\mu t}\bar{x}_0 \\ 0 \end{pmatrix} e^{(B-\lambda I)t}\nabla h(x_0) \\ &\to \begin{pmatrix} 0 & 0 \\ 0 & 1 \end{pmatrix} \nabla h(x_0) \\ &= \begin{pmatrix} \underline{0} \\ D_0 \end{pmatrix} \end{aligned}$$

as $t \to \infty$.

(ii) If $|\phi_t(x)| < \delta/2$, then $|e^{Bt}h(x)| = |h(\phi_t(x))| < \delta$ and so

$$\begin{aligned} |\nabla\phi_t(x)| &= |\nabla h^{-1}(e^{Bt}h(x))e^{Bt}\nabla h(x)| \\ &\leq \sup_{|y|<\delta} |\nabla h^{-1}(y)| \sup_{|y|<\delta} |\nabla h(y)| e^{\lambda t} \\ &< 4e^{\lambda t}. \end{aligned}$$

(iii) Since $h$ and $h^{-1}$ have uniformly Hölder continuous partial derivatives, there exists some $K_0 \in \mathbb{R}$ and $0 < \alpha < 1$ such that

$$|\nabla h(w) - \nabla h(z)| \leq K_0 |w - z|^\alpha$$

and

$$|\nabla h^{-1}(w) - \nabla h^{-1}(z)| \leq K_0 |w - z|^\alpha.$$

Therefore

$$\begin{aligned} &|\nabla\phi_t(x+y) - \nabla\phi_t(x)| \\ &\quad = |\nabla h^{-1}(e^{Bt}h(x+y))e^{Bt}\nabla h(x+y) - \nabla h^{-1}(e^{Bt}h(x))e^{Bt}\nabla h(x)| \\ &\quad \leq |\nabla h^{-1}(e^{Bt}h(x+y))e^{Bt}(\nabla h(x+y) - \nabla h(x))| \\ &\qquad + |(\nabla h^{-1}(e^{Bt}h(x+y)) - \nabla h^{-1}(e^{Bt}h(x)))e^{Bt}\nabla h(x)| \\ &\quad \leq 2e^{\lambda t}K_0|y|^\alpha + 2e^{\lambda t}K_0|e^{\lambda t}(h(x+y) - h(x))|^\alpha \\ &\quad \leq 8K_0 e^{\lambda t(1+\alpha)}|y|^\alpha. \end{aligned}$$ □

Suppose $z \in \mathbb{R}^2$, with $0 < |z| < 1$, and $x_z = x_0 + z$.

THEOREM 3.4. *Fix $C$ and consider the limit $z \to 0$ with $|\frac{z}{D_0 z}| < C$, where $D_0$ is defined in Lemma 3.3. There exist $w_i$, $i = 1, 2$ (not necessarily unique) with $w_i(t, z) \to 0$ uniformly in $t \in [R, -\frac{1}{\lambda}\log|z| - R]$ as $z \to 0$ and $R \to \infty$ such that*

$$\phi_t(x_z) = \bar{x}_0 e^{-\mu t}(e_1 + w_1) + D_0 z e^{\lambda t}(e_2 + w_2).$$

PROOF. Suppose $|\frac{z}{D_0 z}| \le C$ and $R > \frac{1}{\lambda}\log\frac{8}{\delta - 4|\bar{x}_0|}$. If $|x - x_0| \le |z|$ and

$$0 \le t \le \left(\inf_{|x-x_0|\le|z|} \inf\left\{s > 0 : |\phi_s(x)| > \frac{\delta}{2}\right\}\right) \wedge \left(-\frac{1}{\lambda}\log|z| - R\right),$$

then

$$\begin{aligned}
|\phi_t(x)| &\le |\phi_t(x_0)| + |\phi_t(x) - \phi_t(x_0)| \\
&\le 2|\bar{x}_0|e^{-\mu t} + |\nabla\phi_t(x_0 + \theta'(x - x_0))||x - x_0| \\
&\le 2|\bar{x}_0|e^{-\mu t} + 4|z|e^{\lambda t} \\
&\le 2|\bar{x}_0| + 4e^{-\lambda R} \\
&< \frac{\delta}{2},
\end{aligned}$$

where $\theta' \in (0, 1)$. Hence $|\phi_t(x)| < \delta/2$ for all $|x - x_0| \le |z|$ and $t \le -\frac{1}{\lambda}\log|z| - R$. Now

$$\phi_t(x_z) = \phi_t(x_0) + \nabla\phi_t(x_0)z + (\nabla\phi_t(x_0 + \theta z) - \nabla\phi_t(x_0))z$$

for some $\theta \in (0, 1)$ and so, defining

$$w_1(t, z) = \bar{x}_0^{-1}(e^{\mu t}\phi_t(x_0) - \bar{x}_0 e_1)$$

and

$$\begin{aligned}
w_2(t, z) &= (D_0 z)^{-1}(e^{-\lambda t}\nabla\phi_t(x_0)z - D_0 z e_2 + e^{-\lambda t}(\nabla\phi_t(x_0 + \theta z) - \nabla\phi_t(x_0))z), \\
\phi_t(x_z) &= \bar{x}_0 e^{-\mu t}(e_1 + w_1) + D_0 z e^{\lambda t}(e_2 + w_2).
\end{aligned}$$

Then $|w_1| \to 0$ uniformly in $t \ge R$ as $R \to \infty$ by Theorem 3.2, and

$$\begin{aligned}
|w_2| &\le \frac{|z|}{|D_0 z|}\left(\left|e^{-\lambda t}\nabla\phi_t(x_0) - \begin{pmatrix}0\\D_0\end{pmatrix}\right| + Ke^{\lambda\alpha t}|z|^\alpha\right) \\
&\le C\left(\left|e^{-\lambda t}\nabla\phi_t(x_0) - \begin{pmatrix}0\\D_0\end{pmatrix}\right| + Ke^{-\lambda\alpha R}\right) \\
&\to 0
\end{aligned}$$

uniformly in $t \in [R, -\frac{1}{\lambda}\log|z| - R]$ as $R \to \infty$ and $z \to 0$, by Lemma 3.3. □

Since $\phi_t^{-1}(x)$ satisfies (11) with $b$ replaced by $-b$, we may apply Lemma 3.2 and Lemma 3.3 to deduce the existence of $x_\infty$ with $0 < |x_\infty| < \delta/8$ such that $e^{\lambda t}\phi_t^{-1}(x_\infty) \to \binom{0}{\bar{x}_\infty}$ for some $\bar{x}_\infty \in \mathbb{R}$ as $t \to \infty$, and $D_\infty$ such that $e^{-\mu t}\nabla\phi_t^{-1}(x_\infty) \to \binom{D_\infty}{\underline{0}}$ as $t \to \infty$. Suppose that as $z \to 0$, the sign of $D_0 z$ is eventually constant and nonzero. As $\bar{x}_\infty$ has the same sign as $x_{\infty,2}$ (see the proof of Lemma 3.2), we may choose $x_\infty$ such that $\frac{D_0 z}{\bar{x}_\infty} > 0$.

There exists some $t_\infty \geq 0$ such that $\phi_t(x_0)$ does not intersect the line $x_\infty + rD_\infty^*$ for any $t \geq t_\infty$. Let

$$s_z = \inf\{t \geq t_\infty : \phi_t(x_z) = x_\infty + rD_\infty^* \text{ for some } r \in \mathbb{R}\}.$$

THEOREM 3.5. *Fix $C > 0$ and consider the limit $z \to 0$ with $|\frac{z}{D_0 z}| \leq C$. Then*

$$s_z - \frac{1}{\lambda}\log\frac{\bar{x}_\infty}{D_0 z} \to 0$$

*and*

$$\left(\frac{\bar{x}_\infty}{D_0 z}\right)^{\mu/\lambda}(\phi_{s_z}(x_z) - x_\infty) \to \bar{x}_0 \frac{D_\infty^*}{|D_\infty|^2}$$

*as $z \to 0$.*

PROOF. We shall prove this theorem in the case where for $z$ sufficiently small $D_0 z$, $\bar{x}_0 > 0$. The other cases are similar.

Since $\frac{\phi_t(x_0)_2}{\phi_t(x_0)_1} = \frac{e^{\mu t}\phi_t(x_0)_2}{e^{\mu t}\phi_t(x_0)_1} \to \frac{0}{\bar{x}_0} = 0$ as $t \to \infty$, there exists some $T \geq 0$ such that $|\frac{\phi_t(x_0)_2}{\phi_t(x_0)_1}| < 1$ for all $t \geq T$. Let

$$t_z = \inf\{t \geq T : |\phi_t(x_z)_1| = |\phi_t(x_z)_2|\}.$$

By expressing $\phi_t(x_z)$ in the form derived in Theorem 3.4, we may use a similar argument to that in Theorem 2.5 to show

$$t_z - \frac{1}{\lambda + \mu}\log\frac{\bar{x}_0}{D_0 z} \to 0$$

as $z \to 0$. Let $f : B(0,1) \to \mathbb{R}$ be defined by $f(z) = \phi_{t_z}(x_z)_1$. Again as in Theorem 2.5,

$$(D_0 z)^{-\mu/(\lambda+\mu)} f(z) \to \bar{x}_0^{\lambda/(\lambda+\mu)}$$

as $z \to 0$.

Define $g : \mathbb{R}^+ \to \mathbb{R}$ by

$$g(y) = \phi_{t_y'}^{-1}\left(x_\infty + y\frac{D_\infty^*}{|D_\infty|^2}\right)_1,$$

where $t'_y$ is defined in the same way as $t_z$ except for $\phi^{-1}$ instead of $\phi$. [The scaling factor of $|D_\infty|^2$ is chosen so that $D_\infty(y\frac{D^*_\infty}{|D_\infty|^2}) = y$.] Note that $\phi_{s_z}(x_z) = x_\infty + g^{-1}(f(z))\frac{D^*_\infty}{|D_\infty|^2}$.

By a similar argument to above, $y^{-\lambda/(\lambda+\mu)}g(y) \to \bar{x}_\infty^{\mu/(\lambda+\mu)}$ as $y \to 0$. But then

$$\begin{aligned}
&\left|\left(\frac{\bar{x}_\infty}{D_0 z}\right)^{\mu/\lambda} g^{-1}(f(z)) - \bar{x}_0\right| \\
&\quad \le (D_0 z)^{-\mu/\lambda}|\bar{x}_\infty^{\mu/\lambda} g^{-1}(f(z)) - f(z)^{(\lambda+\mu)/\lambda}| \\
&\qquad + |((D_0 z)^{-\mu/(\lambda+\mu)} f(z))^{(\lambda+\mu)/\lambda} - \bar{x}_0| \\
&\quad = \left(\frac{(D_0 z)^{-\mu/(\lambda+\mu)} f(z)}{y^{-\lambda/(\lambda+\mu)} g(y)}\right)^{(\lambda+\mu)/\lambda} |\bar{x}_\infty^{\mu/\lambda} - (y^{-\lambda/(\lambda+\mu)} g(y))^{(\lambda+\mu)/\lambda}| \\
&\qquad + |((D_0 z)^{-\mu/(\lambda+\mu)} f(z))^{(\lambda+\mu)/\lambda} - \bar{x}_0| \\
&\quad \to 0
\end{aligned}$$

as $z \to 0$, where $y = g^{-1}(f(z)) \to 0$ as $z \to 0$. So

$$\left(\frac{\bar{x}_\infty}{D_0 z}\right)^{\mu/\lambda} (\phi_{s_z}(x_z) - x_\infty) = \left(\frac{\bar{x}_\infty}{D_0 z}\right)^{\mu/\lambda} g^{-1}(f(z)) \frac{D^*_\infty}{|D_\infty|^2} \to \bar{x}_0 \frac{D^*_\infty}{|D_\infty|^2}.$$

Also, since $t'_y = s_z - t_z$, and $t'_y - \frac{1}{\lambda+\mu}\log\frac{\bar{x}_\infty}{y} \to 0$ as $y \to 0$,

$$(s_z - t_z) - \frac{1}{\lambda+\mu}\log\frac{\bar{x}_\infty}{(D_0 z/\bar{x}_\infty)^{\mu/\lambda}\bar{x}_0} \to 0,$$

that is,

$$s_z - \frac{1}{\lambda}\log\frac{\bar{x}_\infty}{D_0 z} \to 0. \qquad \square$$

**4. Convergence of the fluctuations.** Now suppose $X_t^N$ is a pure jump Markov process satisfying all the conditions in Section 2, except with $b(x) = Bx + \tau(x)$, $B$ and $\tau$ defined as in Section 3. In this section we express $\phi_{t-t_0}(X_{t_0}^N)$ in a linear form for large values of $N$ and $t_0$.

Recall from Section 2 that $\gamma_t^N = N^{1/2}(X_t^N - x_t)$ and $\gamma_t^N \Rightarrow \gamma_t$ for each $t$ as $N \to \infty$, where $(\gamma_t)_{t\ge 0}$ is the unique solution to the linear stochastic differential equation (4).

Fix some $t_0 \ge 0$. Then $\phi_{t-t_0}(X_{t_0}^N) = \phi_t(\phi_{t_0}^{-1}(X_{t_0}^N))$ and using the same notation as in Section 2, there exists some $\theta \in (0,1)$ such that

$$\begin{aligned}
\phi_{t_0}^{-1}(X_{t_0}^N) &= \phi_{t_0}^{-1}(x_{t_0}) + N^{-1/2}\nabla\phi_{t_0}^{-1}(x_{t_0})\gamma_{t_0}^N \\
&\qquad + N^{-1/2}(\nabla\phi_{t_0}^{-1}(x_{t_0} + \theta N^{-1/2}\gamma_{t_0}^N) - \nabla\phi_{t_0}^{-1}(x_{t_0}))\gamma_{t_0}^N \\
&= x_0 + N^{-1/2} Z_{t_0}^N,
\end{aligned}$$

where $Z^N_{t_0} \Rightarrow Z_{t_0} = \nabla\phi^{-1}_{t_0}(x_{t_0})\gamma_{t_0}$ as $N \to \infty$. Now

$$\begin{aligned} D_0 Z_{t_0} &= \lim_{t\to\infty} e_2^* e^{-\lambda t}\nabla\phi_t(x_0)\int_0^{t_0}\nabla\phi_s^{-1}(x_s)\sigma(x_s)\,dW_s \\ &= \lim_{t\to\infty} e_2^* e^{-\lambda t}\int_0^{t_0}\nabla\phi_{t-s}(x_s)\sigma(x_s)\,dW_s \end{aligned}$$

and

$$\begin{aligned} &\liminf_{t\to\infty} e^{-2\lambda t}\int_0^\infty |e_2^*\nabla\phi_{t-s}(x_s)\sigma(x_s)|^2\,ds \\ &\qquad\le \liminf_{t\to\infty} e^{-2\lambda t}\int_0^\infty |\nabla\phi_{t-s}(x_s)_2|^2|a(x_s)|\,ds \\ &\qquad\le \liminf_{t\to\infty} e^{-2\lambda t}\int_0^\infty 16|D_s|^2 e^{2\lambda(t-s)}A\,ds \\ &\qquad\le \frac{32A}{\lambda}, \end{aligned}$$

where $A$ is defined in (3) and the modulus of $D_s = \lim_{t\to\infty} e^{-\lambda t}\nabla\phi_t(x_s)_2$ is bounded above by 2, by the same argument used to show existence of $D_0$ in Theorem 3.3. Hence, if we define

$$\begin{aligned} \sigma_\infty^2 &= \int_0^\infty \lim_{t\to\infty} e^{-2\lambda t}\nabla\phi_{t-s}(x_s)_2 a(x_s)\nabla\phi_{t-s}(x_s)_2^*\,ds \\ &= \int_0^\infty e^{-2\lambda s}D_s a(x_s)D_s^*\,ds, \end{aligned}$$

then $D_0 Z_{t_0} \to Z_\infty$ almost surely as $t_0 \to \infty$, where $Z_\infty \sim N(0,\sigma_\infty^2)$.

Choose $x_\infty^+$ and $x_\infty^-$, with $0 < |x_\infty^\pm| < \delta/2$ and $x_{\infty,2}^- < 0 < x_{\infty,2}^+$, such that $\phi_t^{-1}(x_\infty^\pm) \to 0$ as $t \to \infty$. Define a random variable $X_\infty$ on the same sample space as $Z_\infty$ by

$$X_\infty = \begin{cases} x_\infty^+, & \text{if } Z_\infty > 0, \\ 0, & \text{if } Z_\infty = 0, \\ x_\infty^-, & \text{if } Z_\infty < 0, \end{cases}$$

and define $\overline{X}_\infty$ similarly, except replacing $x_\infty^\pm$ by $\bar{x}_\infty^\pm$.

By the Skorohod representation theorem, we may assume we are working in a sample space in which $Z^N_{t_0} \to Z_{t_0}$ almost surely for all $t_0 \in \mathbb{N}$. Without this assumption, analogous results about weak convergence hold, however this assumption simplifies the formulation. Let

$$S_{N,t_0} = \inf\{s > t_\infty : \phi_{s-t_0}(X^N_{t_0}) = X_\infty + rD_\infty^* \text{ for some } r \in \mathbb{R}\} \tag{12}$$

and

$$S_N = \frac{1}{2\lambda}\log N + \frac{1}{\lambda}\log\frac{\overline{X}_\infty}{Z_\infty}, \tag{13}$$

where we interpret $\frac{0}{0} = 1$.

THEOREM 4.1. *Suppose* $\sigma_\infty \neq 0$.

(i) *As* $N \to \infty$ *and then* $t_0 \to \infty$,

$$e^{\mu t}|\phi_{t-t_0}(X_{t_0}^N) - \phi_t(x_0)| \to 0$$

*uniformly in* $t$ *on compacts in probability.*

(ii) *If* $R \leq t \leq \frac{1}{2\lambda}\log N - R$, *then there exist* $\varepsilon_i'(N, t_0, t) \to 0$, *uniformly in* $t$ *in probability as* $R, N \to \infty$ *and then* $t_0 \to \infty$ *such that*

$$\phi_{t-t_0}(X_{t_0}^N) = \bar{x}_0 e^{-\mu t}(e_1 + \varepsilon_1') + N^{-1/2} Z_\infty e^{\lambda t}(e_2 + \varepsilon_2').$$

(iii) *As* $N \to \infty$ *and then* $t_0 \to \infty$, $S_{N,t_0} - S_N \to 0$ *in probability. Furthermore, if* $t = S_{N,t_0} - s$ *for some* $s$, *then*

$$e^{\lambda s}|\phi_{t-t_0}(X_{t_0}^N) - \phi_s^{-1}(X_\infty)| \to 0$$

*uniformly in* $s$ *on compacts, in probability as* $N \to \infty$ *and then* $t_0 \to \infty$.

PROOF. (i) By Theorem 3.3, for some $\theta \in (0, 1)$

$$\begin{aligned} e^{\mu t}|\phi_{t-t_0}(X_{t_0}^N) - \phi_t(x_0)| &= e^{\mu t}|\nabla\phi_t(x_0 + \theta N^{-1/2} Z_{t_0}^N)| N^{-1/2}|Z_{t_0}^N| \\ &\leq 4e^{(\lambda+\mu)t} N^{-1/2}|Z_{t_0}^N| \\ &\to 0 \end{aligned}$$

uniformly in $t$ on compacts in probability.

(ii) Apply Theorem 3.4 with $z = N^{-1/2} Z_{t_0}^N$ and use the fact that $D_0 Z_{t_0}^N \to Z_\infty$ almost surely as $N \to \infty$ and then $t_0 \to \infty$. A potential problem arises when $Z_\infty$ is close to 0, however as it is a Gaussian random variable, the probability of this occurring can be made arbitrarily small.

(iii) The first result follows from Theorem 3.5 by a similar argument to (ii). For the second result apply a similar argument to the proof of (i) to $\phi_t^{-1}$. □

**5. A fluid limit for jump Markov processes.** We now show that for large values of $N$ and $t$, $X_t^N$ is in some sense close to $\phi_{t-t_0}(X_{t_0}^N)$ as $t_0 \to \infty$, and combine this with results from Section 3 to obtain results analogous to those in the linear case in Section 2.

Let $f(t, x) = e^{-Bt}(x - \phi_{t-t_0}(X_{t_0}^N))$. By Itô's formula,

$$f(t, X_t^N) = f(0, X_0^N) + M_t^{B,N} + \int_0^t \left(\frac{\partial f}{\partial t} + Kf\right)(s, X_{s-}^N)\, ds,$$

where

$$\begin{aligned}\frac{\partial f}{\partial t} &= -Be^{-Bt}x - e^{-Bt}\tau(\phi_{t-t_0}(X_{t_0}^N)),\\ Kf(s, X_{s-}^N) &= \int_{\mathbb{R}^2}(f(s, X_{s-}^N + y) - f(s, X_{s-}^N))K^N(X_{s-}^N, dy)\\ &= \int_{\mathbb{R}^2} e^{-Bs}yK^N(X_{s-}^N, dy)\\ &= e^{-Bs}b^N(X_{s-}^N)\end{aligned}$$

and

$$\begin{aligned}M_t^{B,N} &= \int_{(0,t]\times\mathbb{R}^2}(f(s, X_{s-}^N + y) - f(s, X_{s-}^N))(\mu^N - \nu^N)(ds, dy)\\ &= \int_{(0,t]\times\mathbb{R}^2} e^{-Bs}y(\mu^N - \nu^N)(ds, dy).\end{aligned}$$

So if $t \geq t_0$, then

$$\begin{aligned}&e^{-Bt}(X_t^N - \phi_{t-t_0}(X_{t_0}^N))\\ &\quad = M_t^{B,N} - M_{t_0}^{B,N}\\ &\quad\quad + \int_{t_0}^t e^{-Bs}(b^N(X_{s-}^N) - b(X_{s-}^N))\,ds\\ &\quad\quad + \int_{t_0}^t e^{-Bs}(\tau(X_{s-}^N) - \tau(\phi_{s-t_0}(X_{t_0}^N)))\,ds.\end{aligned}\tag{14}$$

Since $\tau \in C^2$, $\nabla\tau$ is Lipschitz continuous on the unit disc with Lipschitz constant denoted by $K_0$. In addition to the restrictions on $\delta$ from Section 3, suppose $\delta < \frac{\lambda\mu}{9K_0(\lambda+\mu)}$.

Theorem 5.1. *For all $\varepsilon > 0$,*

$$\lim_{t_0\to 0}\limsup_{N\to\infty}\mathbb{P}\left(\sup_{t_0\leq t\leq S_{N,t_0}} e^{-\lambda t}|X_t^N - \phi_{t-t_0}(X_{t_0}^N)| > \varepsilon N^{-1/2}\right) = 0.$$

Proof. Let

$$R_{N,t_0} = \inf\{t \geq t_0 : e^{-\lambda t}|X_t^N - \phi_{t-t_0}(X_{t_0}^N)| \geq N^{-1/2}\varepsilon\} \wedge S_{N,t_0}.$$

We shall show that $R_{N,t_0} = S_{N,t_0}$ by bounding the terms on the right-hand side of (14).

Fix $c \geq 0$. Since increasing $\varepsilon$ decreases the above probability, we may assume $0 < \varepsilon < \eta_0 \wedge \frac{\lambda e^{-\lambda c}}{9K_0}$. Suppose $C \geq 4$ and pick $R \geq \frac{1}{\lambda}\log(\frac{18CK_0e^{\lambda c}}{\lambda})$.

Define

$$\Omega^1_{N,t_0} = \Big\{\sup_{t\ge t_0} e^{-\lambda t}|e^{Bt}(M^{B,N}_t - M^{B,N}_{t_0})| < N^{-1/2}\frac{\varepsilon}{3}\Big\},$$

$$\begin{aligned}\Omega^2_{N,t_0,R} = &\Big\{\sup_{0\le t\le R} e^{\mu t}|\phi_{t-t_0}(X^N_{t_0}) - \phi_t(x_0)| < \frac{\delta}{2}\Big\}\\ &\cap \Big\{\sup_{R<t<S_{N,t_0}-R} |\varepsilon'_1(N,t_0,t)| \vee |\varepsilon'_2(N,t_0,t)| < 1\Big\}\\ &\cap \Big\{\sup_{S_{N,t_0}-R\le t\le S_{N,t_0}} e^{\lambda(S_{N,t_0}-t)}|\phi_{t-t_0}(X^N_{t_0}) - \phi^{-1}_{S_{N,t_0}-t}(X_\infty)| < \frac{\delta}{2}\Big\},\end{aligned}$$

where $\varepsilon'_1$ and $\varepsilon'_2$ are defined in Theorem 4.1, and

$$\Omega^3_{N,t_0,c} = \Big\{S_{t_0,N} \le \frac{1}{2\lambda}\log N + c\Big\}.$$

Let $N_0$ be sufficiently large that $\sup_{N\ge N_0} N^{1/2}\|b^N - b\| < \lambda\varepsilon/3$.

On the set $\Omega^1_{N,t_0} \cap \Omega^2_{N,t_0,R} \cap \Omega^3_{N,t_0,c} \cap \{C^{-1} < |Z_\infty| < C\}$ with $N \ge N_0$, if $t_0 \le t < R$, then

$$|\phi_{t-t_0}(X^N_{t_0})| \le \delta e^{-\mu t},$$

if $R \le t \le S_{N,t_0} - R$, then

$$\begin{aligned}|\phi_{t-t_0}(X^N_{t_0})| &\le |\bar{x}_0|e^{-\mu t}(1+|\varepsilon'_1|) + N^{-1/2}|Z_\infty|e^{\lambda t}(1+|\varepsilon'_2|)\\ &\le \frac{\delta}{2}e^{-\mu t} + N^{-1/2}2Ce^{\lambda t},\end{aligned}$$

and if $S_{N,t_0} - R \le t \le S_{N,t_0}$, then

$$|\phi_{t-t_0}(X^N_{t_0})| < \delta e^{-\lambda(S_{N,t_0}-t)}.$$

From (14),

$$\begin{aligned}&e^{-\lambda t}|X^N_t - \phi_{t-t_0}(X^N_{t_0})|\\ &\quad\le e^{-\lambda t}|e^{Bt}(M^{B,N}_t - M^{B,N}_{t_0})| + e^{-\lambda t}\int_{t_0}^t |e^{B(t-s)}||b^N(X^N_{s-}) - b(X^N_{s-})|\,ds\\ &\qquad + e^{-\lambda t}\int_{t_0}^t |e^{B(t-s)}||\tau(X^N_{s-}) - \tau(\phi_{s-t_0}(X^N_{t_0}))|\,ds\\ &\quad\le e^{-\lambda t}|e^{Bt}(M^{B,N}_t - M^{B,N}_{t_0})| + \frac{1}{\lambda}\|b^N - b\|\\ &\qquad + \int_{t_0}^t e^{-\lambda s}|\nabla\tau(\phi_{s-t_0}(X^N_{t_0}) + \theta(X^N_{s-} - \phi_{s-t_0}(X^N_{t_0})))||X^N_{s-} - \phi_{s-t_0}(X^N_{t_0})|\,ds\end{aligned}$$

$$\leq e^{-\lambda t}|e^{Bt}(M_t^{B,N} - M_{t_0}^{B,N})| + \frac{1}{\lambda}\|b^N - b\|$$

$$+ K_0 \int_{t_0}^{t} (|\phi_{s-t_0}(X_{t_0}^N)| + |X_{s-}^N - \phi_{s-t_0}(X_{t_0}^N)|) e^{-\lambda s}|X_{s-}^N - \phi_{s-t_0}(X_{t_0}^N)|\, ds,$$

for some $\theta \in (0,1)$.

Hence, on $\Omega^1_{N,t_0} \cap \Omega^2_{N,t_0,R} \cap \Omega^3_{N,t_0,c} \cap \{C^{-1} < |Z_\infty| < C\}$ with $N \geq N_0$,

$$\sup_{t_0 \leq t \leq R_{N,t_0}} e^{-\lambda t}|X_t^N - \phi_{t-t_0}(X_{t_0}^N)|$$

$$\leq N^{-1/2}\frac{\varepsilon}{3} + N^{-1/2}\frac{\varepsilon}{3}$$

$$+ K_0 \int_{t_0}^{R_{N,t_0}} (|\phi_{s-t_0}(X_{t_0}^N)| + |X_{s-}^N - \phi_{s-t_0}(X_{t_0}^N)|) N^{-1/2}\varepsilon\, ds$$

$$\leq N^{-1/2}\varepsilon\biggl(\frac{2}{3} + K_0\biggl(\int_{t_0}^{R_{t_0,N}} (\delta(e^{-\mu t} + e^{-\lambda(S_{N,t_0}-t)}) + N^{-1/2}\varepsilon e^{\lambda t})\, dt$$

$$+ \int_{t_0}^{S_{N,t_0}-R} N^{-1/2} 2Ce^{\lambda t}\, dt\biggr)\biggr)$$

$$\leq N^{-1/2}\varepsilon\biggl(\frac{2}{3} + K_0\biggl(\frac{\delta(\lambda+\mu)}{\lambda\mu} + \frac{\varepsilon e^{\lambda c}}{\lambda} + \frac{2Ce^{\lambda c}}{\lambda}e^{-\lambda R}\biggr)\biggr)$$

$$< N^{-1/2}\varepsilon.$$

Since $X_t^N$ is right continuous, this means $R_{N,t_0} = S_{N,t_0}$ and so

$$\mathbb{P}\Bigl(\sup_{t_0 \leq t \leq S_{N,t_0}} e^{-\lambda t}|X_t^N - \phi_{t-t_0}(X_{t_0}^N)| > N^{-1/2}\varepsilon\Bigr)$$

$$\leq \mathbb{P}((\Omega^1_{N,t_0})^c) + \mathbb{P}((\Omega^2_{N,t_0,R})^c) + \mathbb{P}((\Omega^3_{N,t_0,c})^c) - \mathbb{P}(|Z_\infty| \notin (C^{-1}, C)).$$

Letting $N, t_0, R, C, c \to \infty$ in that order, and using Lemma 2.1 and Theorem 4.1 gives

$$\lim_{t_0\to\infty} \limsup_{N\to\infty} \mathbb{P}\Bigl(\sup_{t \leq S_{N,t_0}} e^{-\lambda t}|X_t^N - \phi_{t-t_0}(X_{t_0}^N)| > N^{-1/2}\varepsilon\Bigr) = 0. \qquad \square$$

REMARK 5.2. The same idea can be used to obtain convergence results for arbitrary matrices $B$, for example, with eigenvalues having the same sign or in higher dimensions. The rate of convergence and the time up to which convergence is valid will depend on the eigenvalues of $B$ and bounds on $|\phi_t(x)|$.

Combining the above result with Theorem 4.1 we get the following.

Theorem 5.3. (i) *For all* $N \in \mathbb{N}$,

$$N^{1/2}|X_t^N - \phi_t(x_0)|$$

*is bounded uniformly in* $t$ *on compacts, in probability. (This follows directly from the fluid limit theorem and diffusion approximation stated in Section 2.)*

(ii) *Suppose* $R \le t \le \frac{1}{2\lambda}\log N - R$. *Then provided* $\sigma_\infty \neq 0$, *for* $i = 1, 2$ *there exist* $\varepsilon_i(N,t) \to 0$ *uniformly in* $t$ *in probability as* $R, N \to \infty$ *such that*

$$X_t^N = \bar{x}_0 e^{-\mu t}(e_1 + \varepsilon_1) + N^{-1/2} Z_\infty e^{\lambda t}(e_2 + \varepsilon_2),$$

[*cf.* (7)].

(iii) *As* $N \to \infty$

$$X_{S_N - s}^N \to \phi_s^{-1}(X_\infty),$$

*uniformly on compacts in* $s \ge 0$, *in probability.*

Remark 5.4. These results can be reformulated as results about weak convergence which are true, independent of the choice of sample space, in a manner analogous to Theorem 2.3. In particular, for any sequence $t_N \to \infty$ as $N \to \infty$, $Z_\infty^N = N^{1/2} e^{-\lambda t_N} X_{t_N,2}^N \Rightarrow Z_\infty$ and working on a space in which this sequence converges almost surely is sufficient for Theorem 5.3.

**6. Continuous diffusion Markov processes.** Interest in this problem arose through looking at the OK Corral problem. It was therefore natural to prove results for pure jump Markov processes. However the proof of the analogous result in the case of continuous diffusion processes is similar and we give it below. The pure jump and continuous cases can be combined to obtain results for more general Markov processes.

Let $(X_t^N)_{t \ge 0}$ be a sequence of diffusion processes, starting from $x_0$ and taking values in some open subset $S \subset \mathbb{R}^2$, that satisfy the stochastic differential equations

$$dX_t^N = \sigma^N(X_t^N)\, dW_t + b^N(X_t^N)\, dt$$

with $\sigma^N, b^N$ Lipschitz.

Suppose that there exist limit functions $b(x) = Bx + \tau(x)$, with $B$ and $\tau$ as in Section 3 and $\sigma$, bounded, satisfying

(a)

$$\sup_{x \in S} N^{1/2}|b^N(x) - b(x)| \to 0.$$

(b)

$$\sup_{x \in S} |N^{1/2}\sigma^N(x) - \sigma(x)| \to 0.$$

It follows that there exists a constant $A$ such that for all $N$

$$\|\sigma^N\| \le (A/N)^{1/2}. \tag{15}$$

Let $\gamma_t^N = N^{1/2}(X_t^N - x_t)$, where $x_t$ is defined as before. It is straightforward, using Gronwall's lemma, to show that $\gamma_t^N \to \gamma_t$ as $N \to \infty$, where $(\gamma_t)_{t\ge 0}$ is the unique solution to the linear stochastic differential equation

$$d\gamma_t = \sigma(x_t)\,dW_t + \nabla b(x_t)\gamma_t\,dt \tag{16}$$

starting from 0, $W$ a Brownian motion.

Consider for $t \ge t_0$ $f(t,x) = e^{-Bt}(x - \phi_{t-t_0}(X_{t_0}^N))$. By Itô's formula,

$$f(t, X_t^N) = f(t_0, X_{t_0}^N) + M_t^{B,N} - M_{t_0}^{B,N} + \int_{t_0}^t \left(\frac{\partial f}{\partial s}(s, X_s^N) + e^{-Bs} b^N(X_s^N)\right) ds,$$

where

$$\frac{\partial f}{\partial t} = -Be^{-Bt}x - e^{-Bt}\tau(\phi_{t-t_0}(X_{t_0}))$$

and

$$M_t^{B,N} = \int_0^t e^{-Bs}\sigma^N(X_s^N)\,dW_s.$$

So if $t \ge t_0$,

$$\begin{aligned} &e^{-Bt}(X_t^N - \phi_{t-t_0}(X_{t_0}^N)) \\ &\quad = M_t^{B,N} - M_{t_0}^{B,N} \\ &\qquad + \int_{t_0}^t e^{-Bs}(b^N(X_{s-}^N) - b(X_{s-}^N))\,ds \\ &\qquad + \int_{t_0}^t e^{-Bs}(\tau(X_{s-}^N) - \tau(\phi_{s-t_0}(X_{t_0}^N)))\,ds. \end{aligned} \tag{17}$$

By comparing this with (14), it is sufficient to prove an analogous result to Lemma 2.1, for the conclusion of Theorem 5.3 to hold for diffusion processes.

LEMMA 6.1. *There exists some constant $C$ such that*

$$\mathbb{E}\left(\sup_{t\ge t_0} e^{-\lambda t}|e^{Bt}(M_t^{B,N} - M_{t_0}^{B,N})|\right) \le CN^{-1/2}e^{-\lambda t_0}.$$

PROOF. By the product rule,

$$\begin{aligned} e^{(B-\lambda I)t}(M_t^{B,N} - M_{t_0}^{B,N}) &= \int_{t_0}^t (B-\lambda I)e^{(B-\lambda I)s}(M_s^{B,N} - M_{t_0}^{B,N})\,ds \\ &\quad + \int_{t_0}^t e^{-\lambda s}\sigma^N(X_s^N)\,dW_s \end{aligned}$$

and hence

$$\begin{aligned}
&\mathbb{E}\Big(\sup_{t\geq t_0} e^{-\lambda t}|e^{Bt}(M_t^{B,N} - M_{t_0}^{B,N})|\Big) \\
&\quad \leq \mathbb{E}\Big(\sup_{t\geq t_0}\int_{t_0}^t (\lambda+\mu)e^{-(\lambda+\mu)s}|(M_s^{B,N} - M_{t_0}^{B,N})_1|\,ds\Big) \\
&\qquad + \mathbb{E}\Big(\sup_{t\geq t_0}\Big|\int_{t_0}^t e^{-\lambda s}\sigma^N(X_s^N)\,dW_s\Big|\Big) \\
&\quad \leq \int_{t_0}^\infty (\lambda+\mu)e^{-(\lambda+\mu)s}(\mathbb{E}(M_s^{B,N} - M_{t_0}^{B,N})_1^2)^{1/2}\,ds \\
&\qquad + \mathbb{E}\Big(\sup_{t\geq t_0}\Big|\int_{t_0}^t e^{-\lambda s}\sigma^N(X_s^N)\,dW_s\Big|^2\Big)^{1/2}.
\end{aligned}$$

Since

$$\mathbb{E}\int_0^t \|e^{-\lambda s}\sigma^N(X_s^N)\|^2\,ds < \infty$$

for all $t \geq 0$, the process

$$\Big(\int_0^t \int_{\mathbb{R}^2} e^{-\lambda s}\sigma^N(X_s^N)\,dW_s\Big)_{t\geq 0}$$

is a martingale, and hence, by Doob's $L^2$ inequality

$$\begin{aligned}
&\mathbb{E}\Big(\sup_{t\geq t_0}\Big|\int_{t_0}^t e^{-\lambda s}\sigma^N(X_s^N)\,dW_s\Big|^2\Big) \\
&\qquad \leq 4\sup_{t\geq t_0}\mathbb{E}\Big(\Big|\int_{t_0}^t e^{-\lambda s}\sigma^N(X_s^N)\,dW_s\Big|^2\Big).
\end{aligned}$$

Now

$$\begin{aligned}
\mathbb{E}((M_t^{B,N} - M_{t_0}^{B,N})_1^2) &= \mathbb{E}\int_{t_0}^t e^{2\mu s}a^N(X_s^N)_{1,1}\,ds \\
&\leq \mathbb{E}\int_{t_0}^t e^{2\mu s}\frac{A}{N}\,ds \\
&\leq \frac{e^{2\mu t}A}{2\mu N},
\end{aligned}$$

where $A$ is defined in (15). Similarly

$$\mathbb{E}\Big(\Big|\int_{t_0}^t e^{-\lambda s}\sigma^N(X_s^N)\,dW_s\Big|^2\Big) \leq \frac{e^{-2\lambda t_0}A}{2\lambda N}.$$

Hence

$$\mathbb{E}\biggl(\sup_{t\ge t_0} e^{-\lambda t}|e^{Bt}(M_t^{B,N}-M_{t_0}^{B,N})|\biggr)$$
$$\le \int_{t_0}^{\infty}(\lambda+\mu)e^{-\lambda s}\biggl(\frac{A}{2\mu N}\biggr)^{1/2}ds + e^{-\lambda t_0}\biggl(\frac{2A}{\lambda N}\biggr)^{1/2}$$
$$\le \frac{A^{1/2}(\lambda+\mu+2(\lambda\mu)^{1/2})}{\lambda(2\mu)^{1/2}}N^{-1/2}e^{-\lambda t_0}. \qquad \square$$

Define $\sigma_\infty, Z_\infty, X_\infty, \overline{X}_\infty$ as in Section 4 and let

$$S_N = \frac{1}{2\lambda}\log N + \frac{1}{\lambda}\log\frac{\overline{X}_\infty}{Z_\infty}.$$

The following analogous theorem to Theorem 5.3 for diffusion processes holds.

THEOREM 6.2. (i) *For all* $N\in\mathbb{N}$,

$$N^{1/2}|X_t^N - \phi_t(x_0)|$$

*is bounded uniformly in* $t$ *on compacts, in probability.*

(ii) *Suppose* $R\le t\le \frac{1}{2\lambda}\log N - R$. *Then provided* $\sigma_\infty\ne 0$, *for* $i=1,2$ *there exist* $\varepsilon_i(N,t)\to 0$ *uniformly in* $t$ *in probability as* $R,N\to\infty$ *such that*

$$X_t^N = \bar{x}_0 e^{-\mu t}(e_1+\varepsilon_1) + N^{-1/2}Z_\infty e^{\lambda t}(e_2+\varepsilon_2),$$

[*cf.* (7)].

(iii) *As* $N\to\infty$

$$X_{S_N-s}^N \to \phi_s^{-1}(X_\infty),$$

*uniformly on compacts in* $s\ge 0$, *in probability.*

**7. Applications.** Throughout this section we work in a sample space in which $Z_\infty^N\to Z_\infty$ almost surely so that, in particular, the statement of Theorem 5.3 holds.

7.1. *Hitting lines through the origin.* As in the linear case, Theorems 5.3 and 6.2 may be used to study the first time that $X_t^N$ hits $l_\theta$ or $l_{-\theta}$, the straight lines passing through the origin at angles $\theta$ and $-\theta$, where $\theta\in(0,\frac{\pi}{2})$, as $N\to\infty$. We define the time that $X_t^N$ first intersects one of the lines $l_{\pm\theta}$ as in Section 2 by

$$T_\theta^N = \inf\biggl\{t\ge 0: \biggl|\frac{X_{t-,2}^N}{X_{t-,1}^N}\biggr| \le |\tan\theta| \text{ and } \biggl|\frac{X_{t,2}^N}{X_{t,1}^N}\biggr| \ge |\tan\theta|\biggr\}.$$

First note that by Lemma 3.2,

$$\frac{\phi_t(x_0)_2}{\phi_t(x_0)_1} = \frac{e^{\mu t}\phi_t(x_0)_2}{e^{\mu t}\phi_t(x_0)_1} \to \frac{0}{\bar{x}_0} = 0$$

as $t \to \infty$. In particular, since $\tan\theta \neq 0$, there exists some $s_\theta \geq 0$ such that $|\frac{\phi_t(x_0)_2}{\phi_t(x_0)_1}| < |\tan\theta|$ for all $t \geq s_\theta$. To rule out the trivial case where $T_\theta^N$ converges to the first time that $\phi_t(x_0)$ hits $l_{\pm\theta}$, we shall assume that $x_0$ is chosen sufficiently close to the origin that $s_\theta = 0$.

We prove the following result in the case where $X_t^N$ is a pure jump process. The proof for continuous diffusion processes is identical, except uses Theorem 6.2 in place of Theorem 5.3.

THEOREM 7.1. *Under the conditions required for Theorem 5.3*

$$T_\theta^N - t_N \Rightarrow c_\theta$$

*and*

$$N^{\mu/(2(\lambda+\mu))}|X^N_{T_\theta^N}| \Rightarrow |\sec\theta||\tan\theta|^{-\mu/(\lambda+\mu)}|\bar{x}_0|^{\lambda/(\lambda+\mu)}|Z_\infty|^{\mu/(\lambda+\mu)}$$

*as $N \to \infty$, where*

$$t_N = \frac{1}{2(\lambda+\mu)}\log N \quad \text{and} \quad c_\theta = \frac{1}{\lambda+\mu}\log\left|\frac{\bar{x}_0\tan\theta}{Z_\infty}\right|.$$

PROOF. By the fluid limit theorem and diffusion approximation, for any constant $R > 0$,

$$\begin{aligned}\mathbb{P}(T_\theta^N \leq R) &\leq \mathbb{P}\left(\sup_{t\leq R}\left|\frac{X_{t,2}^N}{X_{t,1}^N}\right| \geq |\tan\theta|\right)\\ &= \mathbb{P}\left(\sup_{t\leq R}\left|\frac{\phi_t(x_0)_2 + N^{-1/2}\gamma_{t,2}^N}{\phi_t(x_0)_1 + N^{-1/2}\gamma_{t,1}^N}\right| \geq |\tan\theta|\right)\\ &\to 0\end{aligned}$$

as $N \to \infty$.

By an identical argument to Theorem 2.5

$$\mathbb{P}(R \leq T_\theta^N \leq t_N + c_\theta - \varepsilon) \to 0$$

and

$$\mathbb{P}(t_N + c_\theta + \varepsilon \leq T_\theta^N \leq S_N - R) \to 0$$

as $R, N \to \infty$. The result follows immediately. □

REMARK 7.2. As in the linear case, the sign of $Z_\infty$ determines whether $X_t^N$ hits $l_\theta$ or $l_{-\theta}$ at time $T_\theta^N$. Since $Z_\infty$ is a Gaussian random variable with mean 0, each event occurs with probability $\frac{1}{2}$. Furthermore, provided $x_\infty$ is chosen sufficiently close to the origin that $\phi_t^{-1}(x_\infty)$ does not intersect $l_{\pm\theta}$, if $X_t^N$ hits one of the two lines then the probability of it hitting either line again before $S_N$ converges to 0 as $N \to \infty$.

7.2. *Minimum distance from the origin.* Another application is to investigate the minimum distance from the origin that $X_t^N$ can attain for large values of $N$.

THEOREM 7.3. *Under the conditions required for Theorem 5.3,*

$$N^{\mu/(2(\lambda+\mu))} \inf_{t \le S_N} |X_t^N| \Rightarrow \left(\frac{\mu}{\lambda}\right)^{\lambda/(2(\lambda+\mu))} \left(\frac{\lambda}{\mu}+1\right)^{1/2} |\bar{x}_0|^{\lambda/(\lambda+\mu)} |Z_\infty|^{\mu/(\lambda+\mu)}$$

*as* $N \to \infty$.

PROOF. By the fluid limit theorem and diffusion approximation, for any constant $R > 0$,

$$\inf_{t \le R} N^{\mu/(2(\lambda+\mu))} |X_t^N| \ge \inf_{t \le R} N^{\mu/(2(\lambda+\mu))} (|\phi_t(x_0)| - N^{-1/2} |\gamma_t^N|) \to \infty$$

as $N \to \infty$.

By Theorem 5.3,

$$\begin{aligned} &\inf_{R \le t \le t_N - R} N^{\mu/(2(\lambda+\mu))} |X_t^N| \\ &\quad \ge \inf_{R \le t \le t_N - R} (e^{\mu(t_N - t)} |\bar{x}_0| (1 - |\varepsilon_1|) - e^{\lambda(t - t_N)} |Z_\infty| (1 + |\varepsilon_2|)) \\ &\quad \to \infty \end{aligned}$$

in probability as $R, N \to \infty$.

For each $c \ge 0$ there exists some $\varepsilon \to 0$ in probability such that

$$\inf_{S_N - c \le t \le S_N} N^{\mu/(2(\lambda+\mu))} |X_t^N| \ge \inf_{0 \le s \le c} N^{\mu/(2(\lambda+\mu))} (|\phi_s^{-1}(X_\infty)| - \varepsilon) \to \infty.$$

Also

$$\begin{aligned} &\inf_{t_N + R \le t \le 1/(2\lambda) \log N - R} N^{\mu/(2(\lambda+\mu))} |X_t^N| \\ &\quad \ge \inf_{t_N + R \le t \le 1/(2\lambda) \log N - R} e^{\lambda(t - t_N)} |Z_\infty| (1 - |\varepsilon_2|) - e^{\mu(t_N - t)} |\bar{x}_0| (1 + |\varepsilon_1|) \\ &\quad \to \infty \end{aligned}$$

in probability as $R, N \to \infty$.

Finally if $t = t_N + c$, then

$$N^{\mu/(2(\lambda+\mu))}|X_t^N| = N^{\mu/(2(\lambda+\mu))}((e^{-\mu t}\bar{x}_0(1+\varepsilon_{1,1}) + N^{-1/2}Z_\infty e^{\lambda t}\varepsilon_{2,1})^2$$
$$+ (e^{-\mu t}\bar{x}_0\varepsilon_{1,2} + N^{-1/2}Z_\infty e^{\lambda t}(1+\varepsilon_{2,2}))^2)^{1/2}$$
$$\to ((e^{-\mu c}\bar{x}_0)^2 + (e^{\lambda c}Z_\infty)^2)^{1/2}$$

in probability uniformly in $c$ on compact intervals. The right-hand side is minimized when

$$c = \frac{1}{2(\lambda+\mu)}\log\frac{\mu\bar{x}_0^2}{Z_\infty^2\lambda}.$$

Therefore

$$N^{\mu/(2(\lambda+\mu))}\inf_{t\le S_N}|X_t^N| \Rightarrow \left(\frac{\mu}{\lambda}\right)^{\lambda/(2(\lambda+\mu))}\left(\frac{\lambda}{\mu}+1\right)^{1/2}|\bar{x}_0|^{\lambda/(\lambda+\mu)}|Z_\infty|^{\mu/(\lambda+\mu)}$$

as $N\to\infty$. □

EXAMPLE 7.4. Let $(U_t^N, V_t^N)$ be a $\mathbb{Z}^2$ valued process modelling the sizes of two populations of the same species with $U_0^N = V_0^N = N$. The environment that they occupy is assumed to be closed with the initial population density independent of $N$. Each individual reproduces at rate 1. Additionally, the individuals are in competition with each other, a death occurring due to competition over resources at rate $\alpha$ and due to aggression between the populations at rate $\beta$. Hence the transition rates are

$$(u,v)\to\begin{cases}(u+1,v), & \text{at rate } u,\\ (u-1,v), & \text{at rate } \alpha u(u+v-1)/N + \beta uv/N,\\ (u,v+1), & \text{at rate } v,\\ (u,v-1), & \text{at rate } \alpha v(u+v-1)/N + \beta uv/N.\end{cases}$$

Let $X_t^N = (U_t^N, V_t^N)/N$. This gives a sequence of pure jump Markov processes, starting from $x_0 = (1,1)$, with Lévy kernels

$$K^N(x,dy) = Nx_1\delta_{(1/N,0)} + N(\alpha x_1(x_1+x_2-1/N) + \beta x_1x_2)\delta_{(-1/N,0)}$$
$$+ Nx_2\delta_{(0,1/N)} + N(\alpha x_2(x_1+x_2-1/N) + \beta x_1x_2)\delta_{(0,-1/N)}.$$

If we let

$$K(x,dy) = x_1\delta_{(1,0)} + (\alpha x_1^2 + \beta x_1x_2)\delta_{(-1,0)} + x_2\delta_{(0,1)} + (\alpha x_2^2 + \beta x_1x_2)\delta_{(0,-1)}$$

then for $S = (0,2)^2$ and $\eta_0 = 1$,

$$m(x,\theta) = x_1e^{\theta_1} + (\alpha x_1^2 + \beta x_1x_2)e^{-\theta_1} + x_2e^{\theta_2} + (\alpha x_2^2 + \beta x_1x_2)e^{-\theta_2}$$

satisfies

$$\sup_{x\in S^N}\sup_{|\theta|\le\eta_0}\left|\frac{m^N(x,N\theta)}{N} - m(x,\theta)\right|\to 0$$

as $N \to \infty$. Therefore

$$b(x) = m'(x,0) = \begin{pmatrix} x_1(1-\alpha x_1 - (\alpha+\beta)x_2) \\ x_2(1-\alpha x_2 - (\alpha+\beta)x_1) \end{pmatrix}.$$

The deterministic differential equation

$$\dot{\phi}_t(x) = b(\phi_t(x)), \phi_0(x) = x$$

is a special case of the Lotka–Volterra model for two-species competition. See [1] for a detailed interpretation of the parameters $\alpha$ and $\beta$. Further generalizations are discussed in [3].

It is straightforward to check that $b(x)$ is $C^1$ on $S$ and satisfies

$$\sup_{x \in S^N} N^{1/2} |b^N(x) - b(x)| \to 0$$

as $N \to 0$, and that

$$a(x) = \begin{pmatrix} x_1(1+\alpha x_1 + (\alpha+\beta)x_2) & 0 \\ 0 & x_2(1+\alpha x_2 + (\alpha+\beta)x_1) \end{pmatrix}$$

is Lipschitz on $S$.

Now $b(x)$ has a saddle fixed point at $(\frac{1}{2\alpha+\beta}, \frac{1}{2\alpha+\beta})$ and by symmetry any point $x$ on the line $x_1 = x_2$ satisfies $\phi_t(x) \to (\frac{1}{2\alpha+\beta}, \frac{1}{2\alpha+\beta})$ as $t \to \infty$. So under an appropriate translation and rotation, the conditions required for Theorem 5.3 are satisfied, with $\lambda = 1$ and $\mu = \frac{\beta}{2\alpha+\beta}$. [Note that $\sigma_\infty^2 > 0$ since $a(x)$ is positive definite on $S$]. Hence for times $t$ satisfying $t \ll t_N$, where $t_N = \frac{2\alpha+\beta}{4(\alpha+\beta)} \log N$, the two populations coexist with the sizes of both being equal. However at time $t_N + O(1)$ the deterministic approximation breaks down and one side begins to dominate. The previous results give a quantitative description of the behavior of the processes in this region, however we do not go into this here. At time $t = S_N + s = \frac{1}{2} \log N + O(1)$, $X_t^N \to \phi_s^{-1}(X_\infty)$ in probability as $N \to \infty$, where $S_N$ is defined in Theorem 5.1 and $X_\infty$ is defined in Section 4. Now $b(x)$ has stable fixed points at $(\alpha^{-1}, 0)$ and $(0, \alpha^{-1})$ and hence $\phi_s^{-1}(X_\infty)$ converges to one of these two fixed points as $s \to \infty$. For any $\varepsilon \in (0,1)$ we say that a population is $\varepsilon$-extinct if the proportion of the original population that remains is less than $\varepsilon$. Thus for arbitrarily small $\varepsilon$, one of the populations will become $\varepsilon$-extinct at time $\frac{1}{2} \log N + O(1)$.

**Acknowledgments.** Special thanks to James Norris of the Statistical Laboratory at Cambridge University for suggesting the problem and for many helpful discussions. I would also like to thank the anonymous referee who suggested the use of Hartman's theorem in Section 3, which considerably shortened the proofs.

Statistical Laboratory
Centre for Mathematical Sciences
University of Cambridge
Wilberforce Road
Cambridge CB3 0WB
United Kingdom
E-mail: A.G.Turner@statslab.cam.ac.uk